\documentclass[preprint]{elsarticle} 
\usepackage{amsmath,amssymb,amscd}

\newtheorem{theorem}{Theorem}[section]
\newtheorem{proposition}[theorem]{Proposition}
\newtheorem{lemma}[theorem]{Lemma}
\newtheorem{corollary}[theorem]{Corollary}
\newtheorem{remark}[theorem]{Remark}

\newtheorem{definition}[theorem]{Definition}
\newtheorem{conjecture}[theorem]{Conjecture}
\newtheorem{example}[theorem]{Example}

\newcommand{\GS}{\ensuremath{\mathfrak S}}

\newcommand{\la}{\ensuremath{\lambda}}
\newcommand{\N}{\ensuremath{\mathbb N}}
\newcommand{\A}{\ensuremath{\mathbb A}}

\newcommand{\Q}{\ensuremath{\mathbb Q}}

\newcommand{\R}{\ensuremath{\mathbb R}}
\newcommand{\C}{\ensuremath{\mathbb C}}
\newcommand{\pri}{\ensuremath{\smallsetminus}}

\newcommand{\p}{\ensuremath{{\partial}}}

\newcommand{\s}{\sigma}

\newcommand{\al}{\alpha}
\newcommand{\CX}{\C[X]}

\begin{document}

\title{The Strong Factorial Conjecture}
\author{Eric Edo}
\address{ERIM,
University of New Caledonia,
BP R4 - 98851,
Noum\'ea CEDEX,
New Caledonia}
\ead{edo@univ-nc.nc}

\author{Arno van den Essen}
\ead{essen@math.ru.nl}
\address{Faculty  of Science, Mathematics and Computer Science,
Radboud University Nijmegen,
Postbus 9010, 6500 GL Nijmegen,
The Netherlands}

\begin{abstract} In this paper, we present an unexpected link between
the Factorial Conjecture (\cite{EWZ}) and Furter's Rigidity Conjecture (\cite{F4}).
The Factorial Conjecture in dimension $m$ asserts that if a polynomial $f$ in
$m$ variables $X_i$ over $\C$ is such that ${\cal L}(f^k)=0$ for all $k\geq
1$, then $f=0$, where ${\cal L}$ is the $\C$-linear map from
$\C[X_1,\cdots,X_m]$ to $\C$ defined by ${\cal L}(X_1^{l_1}\cdots
X_m^{l_m})=l_1!\cdots l_m!$. The Rigidity Conjecture
asserts that a univariate polynomial map $a(X)$ with complex coefficients
of degree at most $m+1$ such that $a(X)\equiv X$ mod $X^2$, is equal to $X$ if $m$
consecutive coefficients of the formal inverse (for the composition) of $a(X)$ are zero.
\end{abstract}

\maketitle

\section{Presentation}

In Section~2, we recall the Factorial Conjecture from \cite{EWZ}. 
We give a natural stronger version of this conjecture which gives the title of this paper. We also recall the Rigidity Conjecture
from \cite{F4}. 
We present an additive and a multiplicative inversion formula.
We use the multiplicative one to prove that the Rigidity Conjecture is a very particular case of the Strong Factorial Conjecture
(see Theorem~\ref{bridge}). As an easy corollary we obtain a new case of the Factorial Conjecture (see Corollary~\ref{cor:bridge}).
In section~3, we study the Strong Factorial Conjecture in dimension~2. We give a new proof of the Rigidity Conjecture $R(2)$ (see Subsection~3.1)
 using the Zeilberger Algorithm (see \cite{PWZ}). We study the case of two monomials (see Subsection~3.2).
In Section~4 (resp. 5) we shortly give some historical details about the origin of the Factorial Conjecture (resp. the Rigidity Conjecture).

\section{The bridge}

In this section, we fix a positive integer $m\in\N_+$. By $\C^{[m]}=\C[X_1,\ldots,X_m]$,
we denote the $\C$-algebra of polynomials in $m$ variables over $\C$.

\subsection{The Strong Factorial Conjecture}

We recall the definition of the \textit{factorial map} (see \cite{EWZ} Definition~1.2):

\begin{definition}{\rm We denote by ${\cal L}:\C^{[m]}\to\C$ the linear map defined by}
$${\cal L}(X_1^{l_1}\ldots X_m^{l_m})=l_1!\ldots l_m!{\rm\ \ for\ all\ \ } l_1,\ldots l_m\in\N.$$
\end{definition}

\begin{remark}\label{rem:Lperm}
{\rm Let $\s\in\GS_m$ be a permutation of the set $\{X_1,\ldots ,X_m\}$. If we extend $\s$ to an automorphism $\tilde{\s}$
of the $\C$-algebra $\C^{[m]}$, then for all polynomials $f\in\C^{[m]}$, we have ${\cal L}(\tilde{\s}(f))={\cal L}(f)$.}
\end{remark}

\begin{remark}\label{rem:Lmult}
{\rm The linear map ${\cal L}$ is not compatible with the multiplication. Nevertheless, ${\cal L}(fg)={\cal L}(f){\cal L}(g)$
if $f,g\in\C^{[m]}$ are two polynomials
such that there exists an $I\subset\{1,\ldots,m\}$ such that $f\in\C[X_i\,;\,i\in I]$ and $g\in\C[X_i\,;\,i\not\in I]$.}
\end{remark}

We recall the Factorial Conjecture (see \cite{EWZ} Conjecture~4.2).

\begin{conjecture}[Factorial Conjecture $FC(m)$]
For all $f\in\C^{[m]}$,
$$(\forall k\in\N_+)\,{\cal L}(f^k)=0\,\Rightarrow\,f=0.$$
\end{conjecture}

To state some partial results about this conjecture it is convenient to introduce the following notation:

\begin{definition}{\rm We define the \textit{factorial set} as the following subset of $\C^{[m]}$:}
$${\rm F}^{[m]}=\{f\in\C^{[m]}\pri\{0\}\,;\,(\exists k\in\N_+)\,{\cal L}(f^k)\ne 0\}\cup\{0\}.$$
\end{definition}

\begin{remark}
{\rm Let $f\in\C^{[m]}$ be a polynomial, we have $f\in F^{[m]}$ if and only if:
$$(\forall k\in\N_+)\,{\cal L}(f^k)=0\Rightarrow f=0.$$
In other words, the factorial set $F^{[m]}$ is the set of all polynomials satisfying the Factorial Conjecture $FC(m)$
and this conjecture is equivalent to ${\rm F}^{[m]}=\C^{[m]}$.}
\end{remark}

To give a stronger version of this conjecture we introduce the following subsets of $\C^{[m]}$:

\begin{definition}{\rm For all $n\in\N_+$, we consider the following subset of $\C^{[m]}$:
$${\rm F}_n^{[m]}=\{f\in\C^{[m]}\pri\{0\}\,;\,(\exists k\in\{n,\ldots,n+{\cal N}(f)-1\})\,{\cal L}(f^k)\ne 0\}\cup\{0\}$$
where ${\cal N}(f)$ denotes the number of (nonzero) monomials in $f$. We define  \textit{the strong factorial set} as:}
$$F_\cap^{[m]}=\bigcap_{n\in\N_+}F_n^{[m]}.$$
\end{definition}

Since, for all $n\in\N_+$, it's clear that ${\rm F}_n^{[m]}\subset {\rm F}^{[m]}$, the following
conjecture is stronger than the Factorial Conjecture.\\

\begin{conjecture}[Strong Factorial Conjecture $SFC(m)$]
${\rm F}_\cap^{[m]}=\C^{[m]}$.\\
In other words, all polynomials are in the strong factorial set.
\end{conjecture}

\begin{remark}\label{rem:sfc}{\rm Let $n\in\N_+$ be a positive integer.\\
a) Let $f\in\C^{[m]}$ be a polynomial,
$f\in F_n^{[m]}$ if and only if:
$$(\forall k\in\{n,\ldots,n+{\cal N}(f)-1\})\,{\cal L}(f^k)=0 \Rightarrow f=0.$$
b) The equality ${\rm F}_n^{[m]}=\C^{[m]}$ assert that, for all set ${\cal M}=\{M_1,\ldots,M_N\}$ of $N\in\N_+$
distinct unitary monomials, the map $\phi_{\cal M}:\C^N\to\C^N$ defined by $\phi_{\cal M}(x_1,\ldots,x_N)=({\cal L}(f^n),\ldots,{\cal L}(f^{n+N-1}))$
where $f=x_1M_1+\cdots+x_NM_N$ is such that $\phi_{\cal M}^{-1}(0)=\{0\}$, in other words $({\cal L}(f^n),\ldots,{\cal L}(f^{n+N-1}))$ is a regular
system of parameters (see~\cite{F4} \S 2.2.)}
\end{remark}

\begin{remark}\label{rem:mono}
{\rm In the following two cases one can easily prove that $f\in F_\cap^{[m]}$.\\
\hspace*{.3cm}a) ${\cal N}(f)=1$ ($f$ is a monomial).\\
\hspace*{.3cm}b) $f\in\R_{\ge 0}^{[m]}$ (all nonzero coefficients of $f$ are real and positive).}
\end{remark}

\begin{remark}
{\rm The authors of \cite{EWZ} proved that $f\in F_1^{[m]}$ in the following two cases:\\
\hspace*{.3cm}a) ${\cal N}(f)\le 2$ (see Proposition~4.3 in \cite{EWZ}).\\
\hspace*{.3cm}b) $f$ is a linear form (see Proposition~4.11 and Remark 4.12 in \cite{EWZ}).}
\end{remark}

\begin{example}{\rm We consider the polynomial $f=X_1-X_2\in\C^{[2]}$. We have ${\cal N}(f)=2$, and for all $n\in\N_+$,
$${\cal L}(f^n)={\cal L}(\sum_{k=0}^n{n\choose k}X_1^{n-k}(-X_2)^k)=\sum_{k=0}^n{n\choose k}(n-k)!k!(-1)^k
=\left\{\begin{array}{c}
\!0\ \text{ if }n\text{ is odd}\\
n! \text{ if }n\text{ is even}
\end{array}\right.$$
Since $n$ or $n+1$ is even, we deduce that $f\in F_n^{[2]}$. Hence $f\in F_{\cap}^{[2]}$.
This example shows that $SFC(2)$ can not be true in positive characteristic.}
\end{example}

\subsection{The Rigidity Conjecture}

The following conjecture is due to Furter (cf. \cite{F4}).

\begin{conjecture}[Rigidity Conjecture $R(m)$]
 Let $a(X)\in\C[X]$ be a polynomial of degree less of equal to $m+1$ such that $a(X)\equiv X$ mod $X^2$.
Assume that $m$ consecutive coefficients of the formal inverse $a^{-1}(X)$ vanish, then $a(X)=X$.
\end{conjecture}

In \cite{F4}, Furter proved $R(1)$ and $R(2)$ (see Subsection~3.1 for an other proof of $R(2)$), he also studied $R(3)$
but the situation is more complicate in this case. It's natural to introduce, for all $n\in\N_+$, the following statements:

\begin{conjecture}[Partial Rigidity Conjecture $R(m)_n$] Let $a(X)\in\C[X]$ be a polynomial of degree less of equal to $m+1$
 such that $a(X)\equiv X$ mod $X^2$. Assume that the coefficients of $X^{n+1},\ldots,X^{n+m}$ of the formal inverse $a^{-1}(X)$ vanish,
then $a(X)=X$.
\end{conjecture}

From the definitions, it's clear that $R(m)\Leftrightarrow(\forall n\in\N_+)\,R(m)_n$.

\begin{proposition}\label{dual}
For all integers $n\ge 2$, $R(m)_n\Leftrightarrow R(n-1)_{m+1}$.
\end{proposition}

One can prove Proposition~\ref{dual}, using the following easy lemma (which can be use to prove $R(m)_1$).

\begin{lemma}\label{invmod}
Let $n\ge 2$ be an integer and let $a(X),b(X)\in\C[[X]]$ be two formal series such that $a(X)\equiv b(X)$ mod $X^2$.
If $a(X)\equiv b(X)$ mod $X^n$ then $a^{-1}(X)\equiv b^{-1}(X)$ mod $X^n$.
\end{lemma}

{\bf Proof of Proposition~\ref{dual}.} We assume $R(m)_n$ and prove $R(n-1)_{m+1}$, the converse is symmetric.
Let $a(X)\in\C[X]$ be such that ${\rm deg}(a(X))\le n$ and $a(X)\equiv X$ mod $X^2$. Assume that the coefficients of $X^{m+2},\ldots,X^{m+n}$ of the formal inverse $a^{-1}(X)$ vanish. That means there exists $b(X)\in\C[X]$ such that ${\rm deg}(b(X))\le m+1$ and $a^{-1}(X)\equiv b(X)$ mod $X^{m+n+1}$.
By Lemma~\ref{invmod}, we deduce $b^{-1}(X)\equiv a(X)$ mod $X^{m+n+1}$. Applying $R(m)_n$ to the polynomial $b(X)$, we conclude $b(X)=X$
and $a(X)=X$.

\begin{lemma}\label{fid2}
Let $\la\in\C$ be a complex number. The formal inverse of the polynomial $X-\la X^2$ is $X\sum_{k=0}^{\infty}c_k(\la X)^k$ where
$(c_k)_{k\in\N}$ is the sequence of positive integers defined by $c_0=1$ and $c_{k+1}=\sum_{i+j=k}c_ic_j$ for all $k\in\N$.
\end{lemma}

With Lemma~\ref{invmod} and Lemma~\ref{fid2}, it's easy to prove $R(m)_2$ (which is equivalent to $R(1)$ by Proposition~\ref{dual}).\\

The Rigidity Conjecture can be interpreted as an invertibility criterion.
The well-known invertibility criterion for polynomials in one variable is:

\begin{proposition}\label{invcriterion1}
 Let $A$ be $\Q$-algebra and let $a(X)=X(1+a_1X+\cdots +a_mX^m)$ be a polynomial
(where $a_1,\ldots,a_m\in A$).
Then the following statements are equivalent:\\
i)  The polynomial $a(X)$ is invertible for the composition.\\
ii) For all $i\in\{1,\ldots, m\}$, the coefficient $a_i$ is a nilpotent element in $A$.\\
iii) The derivative $a'(X)$ is invertible for the multiplication.
\end{proposition}

\begin{proposition}\label{invcriterion2} If we assume $R(m)$ then, in Proposition~\ref{invcriterion1} we can add:\\
iv)  $m$ consecutive coefficients of the formal inverse $a^{-1}(X)$ are nilpotent.
\end{proposition}

{\bf Proof.} We assume $R(m)$. It's clear that {\it i)} implies {\it iv)}. We prove that {\it iv)} implies {\it ii)}.
Let $a(X)$ be as in Proposition~\ref{invcriterion1}. We assume {\it iv)}.\\
1) First we prove the statement in case $A$ is a domain. We consider the field
$\Q(a_1,\ldots,a_m)$. By Lefschetz principle we can embed this field into $\C$.
Hence we may assume that $a(X)\in\CX$. Since $0$ is the only nilpotent element in $A$, {\it iv)} and $R(m)$ imply $a(X)=X$.\\
2) To prove the general case, let $\wp$ be a prime ideal in $A$. The ring
$\overline{A}=A/\wp$ is a domain.
Since $m$ consecutive coefficients of the formal inverse $a^{-1}(X)$ are nilpotent, the same
holds for $\overline{a^{-1}}(X)$.
It then follows from 1) that $\overline{a}(X)=X$. Hence $\overline{a_i}=0$ for
all $i\in\{1,\ldots,m\}$, {\it i.e.} $a_i\in \wp$
for all $i\in\{1,\ldots,m\}$. Since this holds for every prime ideal $\wp$ of $A$, we obtain that
each $a_i$ belongs to $\cap\wp$, {\it i.e.} each $a_i$ is nilpotent.

\subsection{Inversion formulas}

To study the Rigidity Conjecture it is natural to use inversion formulas to obtain explicit expressions of the coefficients
of the formal inverse. The following formulas are consequences of the Lagrange inversion formula (cf. \cite{St} Corollary~5.4.3).\\

{\bf Lagrange Inversion Formula.}{\it~Let $a(X)\in\C[[X]]$ be a power series such that $a(X)\equiv X$ mod $X^2$.
The formal inverse of $a(X)$ is given by the formula $a^{-1}(X)=X(1+{1\over n+1}\sum_{n\ge 1}u_nX^n)$ where $u_n=[X^n](X^{-1}a(X))^{-(n+1)}$
(the notation $[X^n]$ means "the coefficient of $X^n$ in").}

\begin{lemma}[Additive Inversion Formula]\label{AIF} Let $\al_1,\ldots,\al_m\in\C$ be $m$ complex numbers.
The formal inverse of $a(X)=X(1-(\al_1 X+\cdots+\al_mX^m))$ is given by the following formula
$a^{-1}(X)=X(1+{1\over n+1}\sum_{n\ge 1}u_nX^n)$
where:
$$u_n={1\over n!}\sum_{j_1+2j_2+\cdots+mj_m=n}{(n+j_1+\cdots+j_m)!\over j_1!\ldots j_m!}\al_1^{j_1}\ldots\al_m^{j_m}\hspace{.3cm}(AIF).$$
Moreover if $m\ge 2$ and $\al_1\ne 0$, then ($j_1=n-2j_2-\ldots-mj_m$):
$$u_n={\al_1^n\over n!}\sum_{2j_2+\cdots+mj_m\le n}{(2n-j_2-\ldots-(m-1)j_m)!\over (n-2j_2-\ldots-mj_m)!j_2!\ldots j_m!} x_2^{j_2}\ldots x_m^{j_m}.$$
where $x_k={\al_k/\al_1^k}$ for all $k\in\{2,\ldots,m\}$.\\
In particular, when $m=2$ (with $x=x_2=\al_2/\al_1^2$ and $k=j_2$):
$$u_n={\al_1^n\over n!}\sum_{k\le {n\over 2}}{(2n-k)!\over (n-2k)!k!} x^k.$$
\end{lemma}

{\bf Proof.} By Lagrange Inversion Formula:
$$u_n=[X^n](X^{-1}a(X))^{-(n+1)}=[X^n](1-(\al_1 X+\cdots+\al_mX^m))^{-(n+1)}$$
$$=[X^n]\sum_{0\le j}{n+j\choose n}(\al_1 X+\cdots+\al_mX^m)^j$$
$$=[X^n]\sum_{0\le j}{n+j\choose n}\sum_{j_1+\cdots+j_m=j}{j!\over j_1!\ldots j_m!}(\al_1 X)^{j_1}\ldots(\al_mX^m)^{j_m}$$
$$={1\over n!}[X^n]\sum_{0\le j}(n+j)!\sum_{j_1+\cdots+j_m=j}{1\over j_1!\ldots j_m!}\al_1^{j_1}\ldots\al_m^{j_m}X^{j_1+\cdots+mj_m}$$
$$={1\over n!}\sum_{j_1+2j_2+\cdots+mj_m=n}{(n+j_1+\cdots+j_m)!\over j_1!\ldots j_m!}\al_1^{j_1}\ldots\al_m^{j_m}.$$
\ \\
Using the same kind of computation we obtain (see also \cite{F4} Lemma~1.3):

\begin{lemma}[Multiplicative Inversion Formula]\label{MIF} Let $\mu_1,\ldots,\mu_m\in\C$ be $m$ complex numbers.
The formal inverse of $a(X)=X(1-\mu_1 X)\ldots (1-\mu_m X)$ is given by the following formula
$a^{-1}(X)=X(1+{1\over n+1}\sum_{n\ge 1}u_nX^n)$
where
$$u_n={1\over (n!)^m}\sum_{j_1+\cdots+j_m=n}{(n+j_1)!\ldots(n+j_m)!\over j_1!\ldots j_m!}\mu_1^{j_1}\ldots\mu_m^{j_m}\hspace{.3cm}(MIF).$$
Moreover if $m\ge 2$ and $\mu_1\ne 0$, then ($j_1=n-j_2-\ldots-j_m$):
$$u_n={\mu_1^{n}\over (n!)^m}\sum_{j_2+\cdots+j_m\le n}{(2n-j_2-\ldots-j_m)!(n+j_2)!\ldots(n+j_m)!\over (n-j_2-\ldots-j_m)!\ldots j_m!}y_2^{j_2}\ldots y_m^{j_m}.$$
where $y_k={\mu_k/\mu_1}$ for all $k\in\{2,\ldots,m\}$.\\
In particular, when $m=2$ (with $y=y_2$ and $k=j_2$):
$$u_n={\mu_1^n\over (n!)^2}\sum_{k\le n}{(2n-k)!(n+k)!\over (n-k)!k!} y^k.$$
\end{lemma}

\begin{remark} {\rm\ \\
1) In the case $m=1$, Lemma~\ref{AIF} (or Lemma~\ref{MIF}) implies that the sequence $(c_k)_{k\in\N}$ defined in Lemma~\ref{fid2} is in fact
given by the formula $c_k={1\over k+1}{2k\choose k}$. For all prime numbers $p$ there exists a $k\in\N$ such that $c_k=0$ mod $p$.
Lemma~\ref{fid2} is a formal result and it holds in any commutative ring.
This two facts imply that $R(1)$ does not hold for fields of characteristic $p\ne 0$.\\
2) Since the additive formula contains less terms, it is easier to use to study particular cases like $R(2)$ (see Subsection~3.1).\\
3) Only the multiplicative formula gives a bridge to the Factorial Conjecture (see Subsection~2.4).}
\end{remark}

\subsection{The bridge}

Using Lemma~\ref{AIF} and~\ref{MIF}, we can reformulate $R(m)_n$ in the following way.

\begin{proposition}\label{sequence} Let $\al_1,\ldots,\al_m\in\C$ (resp. $\mu_1,\ldots,\mu_m\in\C$) be $m$ complex numbers
and let $(u_n)_{n\ge 1}$ be the sequence defined by $(AIF)$ (resp. $(MIF)$) in Lemma~\ref{AIF} (resp. Lemma~\ref{MIF}).
For all $n\in\N_+$, the Rigidity Conjecture $R(m)_n$ is equivalent to the following implication:
If $u_n=\ldots=u_{n+m-1}=0$ then $\al_1=\ldots=\al_m=0$ (resp. $\mu_1=\ldots=\mu_m=0$)
\end{proposition}

To state the bridge between the Strong Factorial Conjecture and the Rigidity Conjecture we make one more definition.

\begin{definition}{\rm We consider the following two subsets of $\C^{[m]}$:}
$$E^{[m]}=\{X_1\ldots X_m(\mu_1X_1+\cdots+\mu_mX_m)\,;\,\mu_1,\ldots,\mu_m\in\C\}$$
$$\text{and}\hspace{.4cm}E^{[m]*}=\{f\in E^{[m]}\,;\,{\cal N}(f)=m\}.$$
\end{definition}

Now, we can state the central result of this paper:

\begin{theorem}[the bridge]\label{bridge} Let $n\in\N_+$ be a positive integer,\\
1) the inclusion $E^{[m]}\subset F^{[m]}_n$ implies $R(m)_n$,\\
2) the conjecture $R(m)_n$ implies $E^{[m]*}\subset F^{[m]}_n$,\\
3) the inclusion $E^{[m]}\subset F^{[m]}_n$ is equivalent to $(\forall m'\le m)\,E^{[m']*}\subset F^{[m']}_n$,\\
4) the inclusion $E^{[m]}\subset F^{[m]}_n$ is equivalent to $(\forall m'\le m)\,R(m')_n$.\\
5) The inclusion $E^{[m]}\subset F^{[m]}_\cap$ is equivalent to $(\forall m'\le m)\,R(m')$.
\end{theorem}

The last point of Theorem~\ref{bridge} says that the Rigidity Conjecture $R(k)$ holds for $k=1,\ldots,m$
if and only if the Strong Factorial Conjecture holds for all polynomials of the form
$X_1\ldots X_m(\mu_1X_1+\cdots+\mu_mX_m)\in\C^{[m]}$ where $\mu_1,\ldots,\mu_m\in\C$.
To prove Theorem~\ref{bridge}, we need two more lemmas.

\begin{lemma}\label{lem:bridge1}
Let $\mu_1,\ldots,\mu_m\in\C$ be $m$ complex numbers and let $(u_n)_{n\ge 1}$ be the sequence defined by $(MIF)$.
If $f=X_1\ldots X_m(\mu_1X_1+\cdots+\mu_mX_m)$ then,
for all integers $n\in\N_+$, we have: $(n!)^{m+1}u_n={\cal L}(f^n)$.
\end{lemma}

{\bf Proof.} On one hand:
$$f^n=\sum_{j_1+\cdots+j_m=n}{n!\over j_1!\ldots j_m!}\,\mu_1^{j_1}\ldots\mu_m^{j_m}X_1^{n+j_1}\ldots X_m^{n+j_m}$$
and on the other hand:
$$(n!)^{m+1}u_n=\sum_{j_1+\cdots+j_m=n}n!\,{(n+j_1)!\over j_1!}\ldots{(n+j_m)!\over j_m!}\mu_1^{j_1}\ldots\mu_m^{j_m},$$
$$=\sum_{j_1+\cdots+j_m=n}{n!\over j_1!\ldots j_m!}\,\mu_1^{j_1}\ldots\mu_m^{j_m}{\cal L}(X_1^{n+j_1}\ldots X_m^{n+j_m})$$
We conclude $(n!)^{m+1}u_n={\cal L}(f^n)$ by linearity of ${\cal L}$.

\begin{lemma}\label{lem:bridge2}
Let $\mu_1,\ldots,\mu_m\in\C$ be $m$ complex numbers.
We consider the polynomial $f=X_1\ldots X_m(\mu_1X_1+\cdots+\mu_mX_m)\in E^{[m]}$. Let $I$ be the set of integers $i\in\{1,\ldots,m\}$
such that $\mu_i\ne 0$. We set $m'={\cal N}(f)={\rm card}(I)$.
There exists a unique increasing one to one map $\sigma:\{1,\ldots,m'\}\to I$.
We define:
$$\hat{f}=\prod_{i=1}^{m'}X_i\sum_{i=1}^{m'}\mu_{\sigma(i)}X_i\in\C^{[m']}$$
Then, for all integers $n\in\N_+$, $f\in F_n^{[m]}$ is equivalent to $\hat{f}\in F_n^{[m']}$.
\end{lemma}

{\bf Proof.} Let $k\in\N_+$ be a positive integer. We write $f=(\prod_{i\not\in I}X_i)g$ with
$$g=\prod_{i\in I}X_i\sum_{i\in I}\mu_iX_i=\prod_{i=1}^{m'}X_{\sigma(i)}\sum_{i=1}^{m'}\mu_{\sigma(i)}X_{\sigma(i)}=\tilde{\sigma}(\hat{f})$$
where $\tilde{\sigma}:\C[X_i;i\in I]\to \C^{[m']}$ extend $\sigma$ to an isomorphism of $\C$-algebras.\\
Using Remark~\ref{rem:Lmult} and Remark~\ref{rem:Lperm}, we deduce that
$${\cal L}(f^k)={\cal L}(\prod_{i\not\in I}X_i^k){\cal L}(g^k)=(k!)^{m-m'}{\cal L}(\tilde{\sigma}(\hat{f})^k)
=(k!)^{m-m'}{\cal L}(\hat{f}^k).$$
Let $n\in\N_+$ be an integer. Using Remark~\ref{rem:sfc} a), $f\in F_n^{[m]}$ if and only if
$$(\forall k\in\{n,\ldots,n+m'-1\})\,{\cal L}(f^k)=0 \Rightarrow f=0.$$
This is equivalent to
$$(\forall k\in\{n,\ldots,n+m'-1\})\,{\cal L}(\hat{f}^k)=0 \Rightarrow f=0.$$
Since $m'={\cal N}(\hat{f})$, this last assertion is equivalent to $\hat{f}\in F_n^{[m']}$.\\

{\bf Proof (of Theorem~\ref{bridge}).}\\
1) We assume $E^{[m]}\subset F^{[m]}_n$. We prove $R(m)_n$ using Proposition~\ref{sequence}.
Let $\mu_1,\ldots,\mu_m\in\C$ be $m$ complex numbers. Let $(u_n)_{n\ge 1}$ be the sequence defined by $(MIF)$.
We consider $f=X_1\ldots X_m(\mu_1X_1+\cdots+\mu_mX_m)$. Since $f\in E^{[m]}$, we have $f\in F^{[m]}_n$.
If $u_n=\ldots=u_{n+m-1}=0$ then Lemma~\ref{lem:bridge1} implies ${\cal L}(f^k)=0$ for $k\in\{n,\ldots,n+m-1\}$. Since
${\cal N}(f)\le m$, we deduce $f=0$, {\it i.e.} $\mu_1=\ldots=\mu_m=0$.\\
2) We assume $R(m)_n$. If $f\in E^{[m]*}$ then there exist $\mu_1,\ldots,\mu_m\in\C^*$ such that
$f=X_1\ldots X_m(\mu_1X_1+\cdots+\mu_mX_m)$ and ${\cal N}(f)=m$. We prove that $f\in F^{[m]}_n$ by contradiction.
If $f\not\in F^{[m]}_n$ then ${\cal L}(f^k)=0$ for all $k\in\{n,\ldots,n+m-1\}$. By Lemma~\ref{lem:bridge1},
this implies $u_n=\ldots=u_{n+m-1}=0$ where $(u_n)_{n\ge 1}$ is the sequence defined by $(MIF)$.
Using Proposition~\ref{sequence}, we deduce $f=0$ which is impossible.\\
3) Since $(\forall m'\le m)\,E^{[m']*}\subset E^{[m]}$, $E^{[m]}\subset F_n^{[m]} \Rightarrow (\forall m'\le m)\,E^{[m']*}\subset F^{[m']}_n$.
Conversely, we assume that $(\forall m'\le m)\,E^{[m']*}\subset F^{[m']}_n$.
If $f\in E^{[m]}$, we set $m'={\cal N}(f)\le m$. Since $\hat{f}\in E^{[m']*}\subset F^{[m']}_n$,
Lemma~\ref{lem:bridge2} implies $f\in F_n^{[m]}$.\\
4) By definition, $E^{[m]}\subset F^{[m]}_n$ implies $(\forall m'\le m)\,E^{[m']}\subset F^{[m']}_n$.
By 1), this implies $(\forall m'\le m)\,R(m')_n$.
By 2), this implies $(\forall m'\le m)\,E^{[m']*}\subset F^{[m']}_n$.
And by 3) this implies $E^{[m]}\subset F^{[m]}_n$.\\
5) It's a consequence of the point 4 (applied for all $n$).

\begin{corollary}\label{cor:bridge}
We have: $E^{[m]}\subset F^{[m]}_1\subset F^{[m]}$,
in particular the Factorial Conjecture holds for all $f\in E^{[m]}$.
\end{corollary}

Corollary~\ref{cor:bridge} illustrates the fact that we can use (the point~4 of) the bridge
to obtain a new result about the Factorial Conjecture using a trivial case of the Rigidity Conjecture ($R(m)_1$).
Since $R(1)$ and $R(2)$ are true and respectively equivalent to $R(m)_2$ and $R(m)_3$ we have also $E^{[m]}\subset F^{[m]}_2$ and $E^{[m]}\subset F^{[m]}_3$.

\section{The SFC in dimension~2}

\subsection{A new proof of $R(2)$}

In \cite{F4}, Furter uses the multiplicative formula (Lemma~\ref{MIF}) to prove $R(2)$. Using a computer (but he did not say how), he obtained
the following recurrence relation (where $(u_n)_{n\ge 1}$ is defined by $(MIF)$ in the case $m=2$):\\
$n(n-1)(\mu_1-\mu_2)^2u_n+(n-1)(2n-1)(\mu_1+\mu_2)(\mu_1-2\mu_2)(\mu_2-2\mu_1)u_{n-1}$\\
\hspace*{7cm}$-3(3n-4)(3n-2)\mu_1^2\mu_2^2u_{n-2}=0$.\\

In this subsection, we give a proof of $R(2)$ based on the additive formula.\\

Let $n\in\N_+$ be a positive integer. We prove $R(2)_n$ using the additive version of Proposition~\ref{sequence}.
Let $\al_1,\al_2\in\C$ be two complex numbers. We assume $u_n=u_{n+1}=0$ where $(u_n)_{n\ge 1}$ is the sequence defined by $(AIF)$ in the case $m=2$.
We prove $\al_1=0$ by contradiction. If $\al_1\ne 0$ then, by Lemma~\ref{AIF}, $n!u_n=\al_1^nP_n(x)$ where $x=\al_2/\al_1^2$ and
$$P_n(X)=\sum_{k\le {n\over 2}}{(2n-k)!\over (n-2k)!k!} X^k.$$
The assumption $u_n=u_{n+1}=0$ means that $x$ is a common zero of $P_n(X)$ and $P_{n+1}(X)$.
We can find a recurrence relation between $P_n(X)$, $P_{n+1}(X)$ and $P_{n+2}(X)$ using Zeilberger's algorithm (see \cite{PWZ}).
 This algorithm also produces a "certificate": a rational function in $n$ and $k$ which can be used to write down automatically a proof of the recurrence relation. It is implemented in Maple and Mathematica. For example, using Maple, after downloading the package EKHAD (cf. \cite{PWZ programs}),
the command:
\begin{verbatim}
f:=(n,k)->factorial(2*n-k)*X^k/(factorial(n-2*k)*factorial(k));
zeil(f(n,k),k,n,N);
\end{verbatim}
gives the following recurrence relation:
$$-3(3n+4)(3n+2)X^2P_n(X)-(2n+3)(9X+2)P_{n+1}(X)+(4X+1)P_{n+2}(X)=0.$$
 If someone wants to check this by hand (which is useless because
of the certificate given by Zeilberger's algorithm), he or she has to verify the following equality:
$$ -3(3n+1)(3n-1)k(k-1)-9(2n+1)(2n-k+1)(n-2k+3)k$$
$$ -2(2n+1)(n-2k+3)(n-2k+2)(n-2k+1)+4(2n-k+3)(2n-k+2)(2n-k+1)k$$
$$ +(2n-k+1)(n-2k+3)(n-2k+2)(2n-k+2)=0$$
which contains $5$ terms (against $12$ if we use the multiplicative formula).\\
Since $P_l(0)\ne 0$ for all $l\in\N_+$, changing $n$ to $n-1,\ldots,0$ in the recurrence relation,
we prove that $x$ is a zero of $P_{n-1}(X),\ldots,P_0(X)$. But $P_0(X)$ is a nonzero constant and we get a contradiction.\\
Now since $\al_1=0$, for all $k\in\{n,n+1\}$, by Lemma~\ref{AIF},
$$u_k=\sum_{2j=k}{k+j\choose k}\al_2^j=\left\{\begin{array}{c}
\hspace{.4cm}0\hspace{1cm} \text{ if }k\text{ is odd}\\
{3k/2\choose k}\al_2^{k/2} \text{ if }k\text{ is even}
\end{array}\right.$$
Since $n$ or $n+1$ is even, we deduce $\al_2=0$.\\

The polynomial $P_n(X)$ is a particular case ($m=1$) of the polynomial in \cite{Ed} Lemma~4.1.
It's proportional to the following Gauss hypergeometric function: $\,_2F_1({-n+1\over 2},{-n\over 2};-2n;-4X)$.

\subsection{Two monomials in two variables}

The Strong Factorial Conjecture in the case of a polynomial $f$ in two variables composed of two monomials
is a natural generalization of $R(2)$ (which corresponds to the case
$f\in E^{[2]}=\{\mu_1X_1^2X_2+\mu_2X_1X_2^2,\mu_1,\mu_2\in\C\}$). Since we know that $R(1)$ and $R(2)$ are true,
the bridge (cf. Corollary~\ref{cor:bridge}) implies that $E^{[2]}\subset F_\cap^{[2]}$.\\

In this section, we fix ${\bf a}=(a_1,a_2)\in\N^2$ and ${\bf b}=(b_1,b_2)\in\N^2$ two \emph{distinct} vectors in $\N^2$.
We introduce the set:
$$E^{[2]}_{\{{\bf a},{\bf b}\}}=\{\mu_1 X_1^{a_1}X_2^{a_2}+\mu_2 X_1^{b_1}X_2^{b_2},\mu_1,\mu_2\in\C\}\subset\C^{[2]}.$$
For example, $E^{[2]}=E^{[2]}_{\{(2,1),(1,2)\}}$. We consider the following family of polynomials:

$$P_n(X)=P_{{\bf a},{\bf b},n}(X)=\sum_{k=0}^n{(b_1n+c_1k)!(b_2n+c_2k)!\over k!(n-k)!}X^k.$$
where $(c_1,c_2)={\bf c}={\bf a}-{\bf b}\ne 0$ and $n\in\N$.\\

\begin{remark}
{\rm We have: $P_{{\bf a},{\bf b},n}(X)=P_{\bar{\bf a},\bar{\bf b},n}(X)$
where $\bar{\bf a}=(a_2,a_1)$ and $\bar{\bf b}=(b_2,b_1)$.}
\end{remark}

\begin{theorem}\label{thm:poly}
Let $n\in\N_+$ be a positive integer. We have: $E^{[2]}_{\{{\bf a},{\bf b}\}}\subset F_n^{[2]}$ if and only if, the polynomials $P_{{\bf a},{\bf b},n}(X)$ and $P_{{\bf a},{\bf b},n+1}(X)$ have no common zero in~$\C$.
\end{theorem}

\begin{lemma}\label{lem:poly}
Let $\mu_1,\mu_2\in\C^*$ be two nonzero complex numbers.
We consider $f=\mu_1 X_1^{a_1}X_2^{a_2}+\mu_2 X_1^{b_1}X_2^{b_2}\in E^{[2]}_{\{{\bf a},{\bf b}\}}$, then
$${\cal L}(f^n)=n!\,\mu_2^n\,P_{{\bf a},{\bf b},n}({\mu_1\over\mu_2}).$$
\end{lemma}

{\bf Proof.} From
$$f^n=n!\sum_{k=0}^n{X_1^{a_1k+b_1(n-k)}X_2^{a_2k+b_2(n-k)}\over k!(n-k)!}\mu_1^k\mu_2^{n-k}$$
we deduce:
$${\cal L}(f^n)=n!\sum_{k=0}^n{(b_1n+(a_1-b_1)k)!(b_2n+(a_2-b_2)k)!\over k!(n-k)!}\mu_1^k\mu_2^{n-k}=n!\mu_2^nP_{{\bf a},{\bf b},n}({\mu_1\over\mu_2}).$$

{\bf Proof (of Theorem~\ref{thm:poly}).}
If $E^{[2]}_{\{{\bf a},{\bf b}\}}\not\subset F_n^{[2]}$, then there exists $f\in E^{[2]}_{\{{\bf a},{\bf b}\}}\pri F_n^{[2]}$.
Since $f\not\in F_n^{[2]}$ we have ${\cal N}(f)=2$ (see Remark~\ref{rem:mono}~a) and there exist $\mu_1,\mu_2\in\C^*$ such that
 $f=\mu_1 X_1^{a_1}X_2^{a_2}+\mu_2 X_1^{b_1}X_2^{b_2}$. Now the assumption $f\not\in F_n^{[2]}$ implies ${\cal L}(f^n)={\cal L}(f^{n+1})=0$.
 Using Lemma~\ref{lem:poly}, we deduce that $\mu_1/\mu_2$ is a common zero of $P_{{\bf a},{\bf b},n}(X)$ and $P_{{\bf a},{\bf b},n+1}(X)$.
 Conversely, if $P_{{\bf a},{\bf b},n}(X)$ and $P_{{\bf a},{\bf b},n+1}(X)$ have a common zero $x\in\C$
 then $x\ne 0$ (since $P_{{\bf a},{\bf b},n}(0)\ne 0$) and the polynomial $f=x X_1^{a_1}X_2^{a_2}+ X_1^{b_1}X_2^{b_2}$
belongs to $E^{[2]}_{\{{\bf a},{\bf b}\}}\pri F_n^{[2]}$ by Lemma~\ref{lem:poly}.\\

Theorem~\ref{thm:poly} implies that the Strong Factorial Conjecture in the case of a polynomial in two variables composed with two monomials
is equivalent to the following conjecture.

\begin{conjecture}[Relatively Prime Conjecture $RPC$]
For all ${\bf a},{\bf b}\in\N^2$ such that ${\bf a}\ne{\bf b}$ and all $n\in\N$,
the polynomials $P_{{\bf a},{\bf b},n}(X)$ and $P_{{\bf a},{\bf b},n+1}(X)$ have no common zero in~$\C$.
\end{conjecture}

\begin{proposition}
For all $n\in\N$, the polynomials $P_{{\bf a},{\bf b},n}(X)$ and $P_{{\bf a},{\bf b},n+1}(X)$ have no common zero in~$\C$ in the following cases:\\
1) ${\bf a}=(a,0)$ (with $a\in\N$) and ${\bf b}=(0,1)$,\\
2) ${\bf a}=(a,0)$ (with $a\in\N$) and ${\bf b}=(a,1)$.
\end{proposition}

{\bf Proof.}\\
1) In this case, ${\bf c}=(a,-1)$ and $\displaystyle P_n(X)=\sum_{k=0}^n{(ak)!\over k!}X^k$. We have the relation:
$\displaystyle P_{n+1}(X)-P_n(X)={(a(n+1))!\over (n+1)!}X^{n+1}$.
If $\alpha\in\C$ is a common root of $P_n$ and $P_{n+1}$ then $\alpha=0$ which is impossible since $P_n(0)=1$.\\
2) In this case, ${\bf c}=(0,-1)$ and $\displaystyle P_n(X)=(an)!\sum_{k=0}^n{1\over k!}X^k$. We have the relation:
$\displaystyle (an)!P_{n+1}(X)-(a(n+1))!P_n(X)={(an)!(a(n+1))!\over (n+1)!}X^{n+1}$.
If $\alpha\in\C$ is a common root of $P_n$ and $P_{n+1}$ then $\alpha=0$ which is impossible since $P_n(0)=(an)!$.\\

For some small values of $a_1,a_2,b_1,b_2$ the polynomial $P_{n}(X)$ satisfies
a recurrence relation of order $2$ given by Zeilberger's algorithm and we deduce
that $P_n(X)$ and $P_{n+1}(X)$ have no common zero in $\C$.
For example, in the case:
${\bf a}=(1,1)$ and ${\bf b}=(0,0)$, we have ${\bf c}=(1,1)$ and $\displaystyle P_n(X)=\sum_{k=0}^n{k!\over (n-k)!}X^k$.
Zeilberger's algorithm gives $XP_n(X)-(n+2)XP_{n+1}(X)+P_{n+2}(X)=0$.
Unfortunately, for bigger values the polynomial $P_n(X)$ satisfies a recurrence relation of order $\ge 3$
and we can't directly deduce that $P_n(X)$ and $P_{n+1}(X)$  have no common zero  in $\C$.
For example, in the case: ${\bf a}=(3,0)$ and ${\bf b}=(0,0)$, we have ${\bf c}=(3,0)$ and
$\displaystyle P_n(X)=\sum_{k=0}^n{(3k)!\over k!(n-k)!}X^k$.
Zeilberger's algorithm gives the relation:
$27XP_n(X)-54(n+2)XP_{n+1}(X)+3(3n+8)(3n+7)XP_{n+2}(X)-P_{n+3}(X)=0$.

\section{The origin of Factorial Conjecture}

In this section, we fix a positive integer $m\in\N_+$.

\subsection{The Jacobian Conjecture}

We recall the famous Jacobian Conjecture proposed by Keller in 1939.

\begin{conjecture}[Jacobian Conjecture $JC(m)$] An endomorphism of $\C^{[m]}$
is an automorphism if and only if the determinant of his Jacobian matrix is a nonzero constant.
\end{conjecture}

We cannot discus here all the details of this fascinating conjecture.
Facing the difficulty of this question, a lot of people discovered different kind of reductions and reformulations
(see for example \cite{BCW}, \cite{Es1} and~\cite{B}). Often, a new conjecture $X_n$ depending on a parameter $n$
is introduced and a result like "$X_m$ is true for all $m$ if and only if $JC(m)$ is true for all $m$"
is obtained (see for example Zhao's Theorem below).

\subsection{The Vanishing Conjecture}

In 2007, Zhao introduced the Vanishing Conjecture. There exists several versions of this conjecture in the literature
(see \cite{Z1} Conjecture 7.1, \cite{Z3} Conjecture 1.1 and \cite{Es2} \S 2). Here is a very particular version:\\

\begin{conjecture}[Vanishing Conjecture $VC(m)$]
Let $\Delta=\p_1^2+\cdots+\p_m^2$ be the Laplace operator.
Let $f\in \C^{[m]}$ be a homogeneous polynomial of degree $4$.
If $\Delta^k (f^k)=0$ for all $k\in\N_+$ then there exists $K\in\N_+$ such that  $\Delta^k (f^{k+1})=0$ for all $k\ge K$.
\end{conjecture}

Zhao proved the equivalence between this conjecture and the Jacobian Conjecture (see \cite{Z1}~Theorem~7.2 and \cite{Z3}~Theorem~1.2).

\begin{theorem}[Zhao]\label{thm:JCVC}
The following two statements are equivalent:\\
(1) For all $m\in\N_+$, the Jacobian Conjecture $JC(m)$ holds.\\
(2) For all $m\in\N_+$, the Vanishing Conjecture $VC(m)$ holds.
\end{theorem}

This theorem is impressive because the Vanishing Conjecture looks very particular
(we have only to deal with homogeneous polynomials of degree $4$) and quite simple (we just have
to check that something is zero). But this kind of result have a weakness because is not "$m$ to $m$" but "for all $m$ to for all~$m$".
For example, if the Jacobian Conjecture turns out to be true (or provable or even less harder to prove) only in dimension $2$, then an equivalence
"for all~$m$" to an other conjecture is completely useless. Nevertheless, the Vanishing Conjecture is interesting by itself
and should be studied even without view to the Jacobian Conjecture.
The Vanishing Conjecture has been generalized in the following way (see \cite{Es2} \S 2):

\begin{conjecture}[Generalized Vanishing Conjecture $GVC(m)$]
Let $\Delta\in \C[\p_1,\ldots,\p_m]$ be a differential operator. Let $f,g\in \C^{[m]}$ be two polynomials.
If $\Delta^k (f^k)=0$ for all $k\in\N_+$ then there exists $K\in\N_+$ such that  $\Delta^k (gf^k)=0$ for all $k\ge K$.
\end{conjecture}

If in this conjecture we replace $\Lambda$ by the Laplace operator $\Delta$ and $f$ by
a homogeneous polynomial of degree $4$, then we obtain $VC(m)$.

\subsection{The Image Conjecture}

In this subsection, we consider  the $\C$-algebra $A=\C[\xi,z]$ where $\xi=(\xi_1,\ldots,\xi_m)$ and $z=(z_1,\ldots,z_m)$ are two sets of commuting indeterminates. Zhao introduced the following concept. A $\C$-linear subspace $M$ of $A$
is called a {\it Mathieu subspace} of $A$ if for every $f$ in $A$ we have:
$$(\forall k\in\N_+)\, f^k\in M \text{ implies } (\forall g\in A)(\exists K\in\N)(\forall k\ge K) gf^k\in M.$$
For all $i\in\{1,\ldots,m\}$, we consider the operator ${\cal D}_i=\xi_i-\p_{z_i}$.
We denote by ${\rm Im}\,{\cal D}$ the image of the map ${\cal D}=({\cal D}_1,\ldots,{\cal D}_m)$ from $A^m$ to $A$ {\it i.\,e.}
$\sum_{i=1}^m{\cal D}_i A$.\\

In 2010, Zhao introduce the following conjecture (see \cite{EWZ} Conjecture~3.1).
We don't discuss here the General Image Conjecture (see \cite{Z3} Conjecture~1.3 or \cite{Es2} \S 1).\\

\begin{conjecture}[Special Image Conjecture $SIC(m)$]
The image ${\rm Im}\,{\cal D}$ is a Mathieu subspace of $A$.
\end{conjecture}

We define the \textit{special image set} as the set of $f\in A$ such that
$$(\forall k\in\N_+)\, f^k\in {\rm Im}\,{\cal D} \text{ implies } (\forall g\in A)(\exists K\in\N)(\forall k\ge K) gf^k\in {\rm Im}\,{\cal D}.$$
Zhao's theorem about the Vanishing Conjecture (see Subsection~4.2) implies that we can add in Theorem~\ref{thm:JCVC} the following statement (see \cite{Z3} Theorem 3.7 and \cite{EWZ} Theorem~3.2):\\

{\it (3) For all  $m\in\N_+$, the polynomials of the form $(\xi_1^2+\cdots+\xi_m^2)P(z)$ with $P(z)\in\C[z]$ homogenous of degree $4$
are in the special image set.}

\subsection{The factorial conjecture}

Let ${\cal E}:\C[\xi,z]\rightarrow\C[z]$ be the $\C$-linear map defined by
$${\cal E}(\xi_1^{i_1}\ldots \xi_m^{i_m}z_1^{j_1}\ldots z_m^{j_m})=\p_{z_1}^{i_1}(z_1^{j_1})\ldots \p_{z_m}^{i_m}(z_m^{j_m})$$

Zhao proved (see \cite{Z3} Theorem~3.1):

\begin{theorem}[Zhao]\label{thm:imker}
${\rm Im}\,{\cal D}={\rm Ker}\,{\cal E}$.
\end{theorem}

Let $f\in\C^{[m]}$ be a polynomial. Since ${\cal E} (\xi_1^{i_1}z_1^{i_1}\ldots\xi_m^{i_m}z_m^{i_m})=i_1!\ldots i_m!$
 \noindent then one readily verifies that ${\cal E}(f(\xi_1z_1,\ldots,\xi_mz_m))={\cal L}(f(X_1,\ldots,X_m))$
 where ${\cal L}$ is the factorial map (see Subsection~2.1). By Theorem~\ref{thm:imker}, the polynomial $f(\xi_1z_1,\ldots,\xi_mz_m)$
is in the special image set if and only if
$$(\forall k\in\N_+)\, {\cal L}(f^k)=0 \text{ implies } (\forall g\in A)(\exists K\in\N)(\forall k\ge K) {\cal L}(gf^k)=0.$$

Since the right part of this implication is automatically true if $f=0$,
the Factorial Conjecture (see Subsection~2.1) implies that all polynomials of the form $f(\xi_1z_1,\ldots,\xi_mz_m)$ with $f\in\C^{[m]}$ are in the special image set.\\

In \cite{EWZ} the authors proved that various polynomials are in the special image set.
They showed that the truth of the one dimensional Image Conjecture is equivalent to the fact that
all polynomials of the form $f(\xi_1z_1)$ are in the special image set.
In dimension greater than $1$ such an equivalence does not hold, but this was the motivation to study the Image Conjecture
in the case of polynomials of the previous form.

\section{The origin of the Rigidity Conjecture}

In this section, we denote by ${\cal G}$ the group of polynomial automorphisms of the complex plane
$\A^2_\C$, by ${\cal A}$ the affine subgroup and by ${\cal B}$ the triangular subgroup.

\subsection{Polydegree and ind-variety topology on $\A_{\C}^2$}

In 1942, Jung  proved that ${\cal G}$ is generated by ${\cal A}$ and ${\cal B}$ (cf. \cite{J}). In
1953, van der Kulk  generalized this result for a field of arbitrary characteristic and
implicitly obtained that ${\cal G}$ is the amalgamated product of ${\cal A}$ and ${\cal B}$ along ${\cal A}\cap{\cal B}$ (cf. \cite{K}).\\

In 1982, Shafarevich endowed the group ${\cal G}$ with the structure of an
infinite-dimensional algebraic variety (cf. \cite{Sh}). If ${\cal H}\subset {\cal G}$, we denote by $\overline{\cal H}$ the
closure of ${\cal H}$ in ${\cal G}$ for the Zariski topology associated with
this structure.\\

In 1989, Friedland and Milnor used the amalgamated product structure of ${\cal G}$
to define the {\it polydegree} of an automorphism $\s\in{\cal G}$
as the sequence of the degrees of the triangular automorphisms in a decomposition of $\s$ as a product of elements
in ${\cal A}\cup{\cal B}$ (cf. \cite{FM}). They proved that the set ${\cal G}_d$ of all automorphisms whose multidegree is $d=(d_1,\ldots,d_l)$
(where $d_1,\ldots,d_l\ge 2$ are integers) is an analytic variety of dimension $d_1+\cdots+d_l+6$.

\subsection{Polydegree Conjectures}

In 1997, Furter separately defined the polydegree (he called it {\it multidegree}) and studied the set ${\cal G}_d$ (cf.~\cite{F1}).
In 2002, he gave the first deep result of this theory: the length of an automorphism is lower semicontinuous (cf. \cite{F2} Theorem~1).
He introduced a partial order $\preceq$ to describe the closure $\overline{{\cal G}_d}$ of the set of automorphism of a fixed polydegree~$d$.
He conjectured that $\overline{{\cal G}_d}$ is a the union of all ${\cal G}_e$ with $e\preceq d$ (cf.~\cite{F2} Section~1~b).\\

In 2004, Furter and the first author used some polynomial automorphisms constructions (developed in~\cite{EV}) to produce counterexample
to this conjecture for polydegree $d$ of length $3$ (cf.~\cite{EF}). Nevertheless, in the length $2$ case it is still open
(see \cite{F3} and also~\cite{Ed}~Conjecture~2):

\begin{conjecture}[The Length $2$ Polydegree Conjecture $PC(m,n)$]
$$\overline{{\cal G}_{(m+1,n+1)}}=\bigcup_{m'\le m\,,\,n'\le n}{\cal G}_{(m'+1,n'+1)}\,\cup\,\bigcup_{k\le m+n+1}{\cal G}_{(k+1)}.$$
\end{conjecture}

In this equality, $m'$ and $n'$ are regarded as integers $\ge 1$ but $k$ starts at $0$ and ${\cal G}_1={\cal A}$ by convention.\\

In 2007, the first author obtained the following partial results of this conjecture (cf.~\cite{Ed}):

\begin{theorem}[Edo]
Let $m,n\in\N_+$ be integers.\\
1) If $n\in m\,\N$, then: ${\cal G}_{(m+n+1)}\subset\overline{{\cal G}_{(m+1,n+1)}}.$\\
2) If $m$ is even and $n\in{m\over 2}\N$, then: ${\cal G}_{(m+n+1)}\cap\overline{{\cal G}_{(m+1,n+1)}}\ne\emptyset.$
\end{theorem}

Recently, Furter introduced the Rigidity Conjecture (cf. \cite{F4}):

\begin{conjecture}[Rigidity Conjecture $R(m,n)$]
Let $a(X),b(X)\in\C[X]$ be two polynomials such that $a(X)\equiv b(X)\equiv X$ mod $X^2$, ${\rm deg}(a)\le m+1$ and
${\rm deg}(b)\le n+1$. If $a\circ b(X)=X$ modulo $X^{m+n+2}$ then $a(X)=b(X)=X$.
\end{conjecture}

The Rigidity Conjecture $R(m)$ is equivalent to $(\forall n\in\N_+)\,R(m,n)$ (cf. \cite{F4}).
Furter proved the following impressive result (cf. \cite{F4}):

\begin{theorem}[Furter]
For all $m,n\in\N_+$, the Rigidity Conjecture $R(m,n)$ implies the Length $2$ Polydegree Conjecture $PC(m,n)$.
\end{theorem}

This gives us a very good reason to study the Rigidity Conjecture!


\begin{thebibliography}{99}


\bibitem{BCW} H. Bass, E. Connell, D. Wright, {\it The Jacobian conjecture: reduction of degree and formal expansion of the inverse.}
 Bull. Amer. Math. Soc. {\bf 7} (1982), no. 2, 287–-330.

\bibitem{B} M. de Bondt, {\it Homogeneous  Keller Maps,}
Ph. D. thesis (2009), Radboud University, Nijmegen, The Netherlands.

\bibitem{Ed} E. Edo, {\it Some families of polynomial automorphisms II},
Acta Math. Vietnam. {\bf 32} (2007), no. 2-3, 155–-168.

\bibitem{EF} E. Edo, J.-P. Furter, \textit{Some families of polynomial automorphisms},
J.~Pure Appl. Algebra {\bf 194} (2004), no. 3, 263--271.

\bibitem{EV} E. Edo, S. V\'{e}n\'{e}reau,
{\it Length 2 variables of} $A[x,y]$ {\it and transfer},
Polynomial automorphisms and related topics (Krak\`{o}w, 1999).
Ann. Polon. Math. {\bf 76} (2001), no. 1-2, 67-–76.

\bibitem{Es1} A. van den Essen, {\it Polynomial automorphisms and the Jacobian conjecture},
Progress in Mathematics, {\bf 190}, Birkh\"auser Verlag, Basel (2000).

\bibitem{Es2} A. van den Essen, {\it The Amazing Image Conjecture},\\
http://arxiv.org/pdf/1006.5801.pdf (2010).


\bibitem{EWZ} A. van den Essen, D. Wright, W. Zhao, {\it On the Image Conjecture},
J.~Algebra {\bf 340} (2011), 211–-224.

\bibitem{FM} S. Friedland, J. Milnor, {\it Dynamical properties of plane polynomial automorphisms},
Ergod. Th. Dyn. Syst {\bf 9} (1989), 67--99.

\bibitem{F1} J.-P. Furter, {\it On the variety of automorphisms of the affine plane},
J.~Algebra {\bf 195} (1997), 604–-623.

\bibitem{F2} J.-P. Furter, {\it On the length of polynomial automorphisms of the affine plane},
Math. Ann. {\bf 322} (2002), no. 2, 401–-411.

\bibitem{F3} J. P. Furter, {\it Plane polynomial automorphisms of fixed multidegree}.
Math. Ann. {\bf 343} (2009), no. 4, 901–-920.

\bibitem{F4} J. P. Furter, {\it Polynomial Composition Rigidity and Plane Polynomial Automorphisms},
http://perso.univ-lr.fr/jpfurter/

\bibitem{J} H.W.E. Jung. {\it \"{U}ber ganze birationale Transformationen der Ebene},
J.~Reine Angew. Math. {\bf 184} (1942), 161--174.

\bibitem{K} W. van der Kulk. {\it On polynomial rings in two variables},
Nieuw. Arch. Wisk. (3) {\bf 1} (1953), 33--41.

\bibitem{PWZ} M. Petkovšek, H. Wilf, D. Zeilberger, {\it A=B},
A K Peters, Ltd., Wellesley, MA (1996). ISBN: 1-56881-063-6.

\bibitem{PWZ programs} http://www.math.rutgers.edu/~zeilberg/programsAB.html

\bibitem{Sh} I.R. Shafarevich. {\it On some infinite-dimensional groups II},
Math. USSR Izv, {\bf 18} (1982), 214--226.

\bibitem{St}  R. Stanley. {\it Enumerative Combinatorics},
Vol. 2, Cambridge Studies Adv. Math. {\bf 62} (1999), Cambridge Univ. Press, Cambridge.

\bibitem{Z1} W. Zhao, {\it Hessian nilpotent polynomials and the Jacobian conjecture},
Trans. Amer. Math. Soc. {\bf 359} (2007), no. 1, 249--274.

\bibitem{Z2} W. Zhao, {\it Generalizations of the image conjecture and the Mathieu conjecture},
J.~Pure Appl. Algebra {\bf 214} (2010), no. 7, 1200–-1216.

\bibitem{Z3} W. Zhao, {\it Images of commuting differential operators of order one with constant leading coefficients},
J.~Algebra {\bf 324} (2010), no. 2, 231--247.

\bibitem{Z4} W. Zhao, {\it Mathieu subspaces of associative algebras},
J.~Algebra {\bf 350} (2012), 245--272.

\end{thebibliography}
\end{document}